\documentclass[12pt]{article}

\usepackage{amssymb,amsmath}
\date{ }
\def\ma{\mathbb}
\def\BE {\begin{eqnarray}}
\def\EE {\end{eqnarray}}
\def\BC {\begin{eqnarray*}}
\def\EC {\end{eqnarray*}}
\def\OPLUS#1{\raisebox{-5pt}{\mbox{$\begin{array}{c}
\oplus\\[-5pt]\scriptstyle#1
\end{array}$}}}

\def\vs{\vspace*}\def\cl{\centerline}
\begin{document}
\cl{{\bf\large Classification of  modules of the intermediate
series}} \cl{\bf\large {over Ramond $N=2$ superconformal
algebras}\footnote{Supported by Communication University of China,
NSF grants 10571120, 10471096 of China and ``One Hundred Talents
Program" from University of Science and Technology of China
\\\indent Correspondence author: J.~Fu: {\it fujy@cuc.edu.cn} }\vs{7pt}}

\cl{\bf Jiayuan Fu} \cl{\small\it Department of Mathematics,
Communication University of China} \cl{\small\it Beijing 100024,
China} \cl{\small\it Email: \vs{7pt}fujy@cuc.edu.cn}

 \cl{\bf Qifen Jiang} \cl{\small\it Department of Mathematics,
Shanghai Jiaotong University, Shanghai 200240, China} \cl{\it
Email:\vs{7pt} qfjiang@sjtu.edu.cn}

  \cl{\bf Yucai Su}
\cl{\small\it Department of Mathematics, University of Science and
Technology of China} \cl{\small\it Anhui 230026, China}\cl{\it
Email: \vs{9pt}ycsu@ustc.edu.cn}

\cl{\bf Abstract}
  \hskip 1pt{\small
   In this paper, we first discuss the structure of the Ramond $N=2$
   superconformal algebras. Then we classify the modules
   of the intermediate series over Ramond $N=2$ superconformal algebra. }
 \vskip 9pt\cl{\bf 1. Introduction}
\def\theequation{1.\arabic{equation}}

More than two decades ago, superconformal algebras were first
constructed independently and almost at the same time by Kac
{\cite{K1}} and by Ademollo et al.~{\cite{A1}}. On the mathematical
side Kac and van de Leuer \cite{KL}, Cheng and Kac \cite{CK}  have
classified all possible superconformal algebras and Kac recently has
proved that their classification is complete.

The Neveu-Schwarz, the Ramond and the Topological $N=2$
superconformal algebras are connected to each other by the spectral
flows and/or the topological twists. Therefore, we only consider the
{\it Ramond $N=2$ superconformal algebra},
 which is a ${\mathbb{Z}}_2$-graded space:
${\cal L}={\cal L}_{\overline{0}}\oplus{\cal L}_{\overline{1}},$
with
$$ {\cal L}_{\overline{0}}=\mbox{span}_{\ma C}\{ L_i, H_j, c\mid i,
j\in{\ma Z}\},  ~~~  {\cal L}_{\overline{1}}=\mbox{span}_{\ma
C}\{G_i^-, G_j^+\mid i, j\in{\ma Z} \},$$ such that $c$ is a central
element and the following relations hold: \BE\label{a23}
\begin{array}{llllll}
 [ L_i, L_j ]=
 (i-j)L_{i+j}+\frac{1}{12}(\!\!\!\!\!&i^3-i)\delta_{i+j,0}c, \\[7pt]
   [ L_i, H_j ]=-jH_{i+j}, &
   [ H_i, H_j ]= \frac{1}{3} i \delta_{i+j,0}c,\\[7pt]
    [ L_i, G_j^{\pm} ]= (\frac{i}{2}-j)G_{i+j}^{\pm},&
   [ H_i, G_j^{\pm} ]=  \pm G_{i+j}^{\pm}, \\[7pt]
   [ G_i^+, G_j^+ ]= [ G_i^-, G_j^- ]= 0,\!\!\!\!\!\!
  &
  [ G_i^-, G_j^+ ]=
  2L_{i+j}-(i-j)H_{i+j}+\frac{1}{3}(i^2-\frac{1}{4})\delta_{i+j,0}c.\!\!\!\!\!\!\!\!\!
 \end{array} \EE
Obviously, the Cartan subalgebra of ${\cal L}$ is ${\cal H}={\mathbb
C}L_0+{\mathbb C}H_0+{\mathbb C}c$, and $Vir=$ span$_{\mathbb C} \{
L_m, c\mid m\in{\mathbb Z}\}$ is a Virasoro subalgebra of ${\cal
L}$, which can be described as the universal central extensions of
the Lie algebras of differential operators (see \cite{A} for
details).

An ${\cal L}$-module $V$ is called a {\it Harish-Chandra module} if
$V$ is a direct sum of its finite dimensional weight spaces
$V^\lambda=\{v\in V\mid x\cdot v=\lambda(x)v, x\in{\cal H} \}$ for
all $\lambda\in{\cal H}^*$ (the dual of ${\cal H}$). Similar to the
case of Virasoro algebra, we can define the module of the
intermediate series over ${\cal L}$:
\\ [5pt]
{\bf Definition 1.1}\quad A {\it module of the intermediate series}
over ${\cal L}$ is an indecomposable Harish-Chandra module $V$ such
that dim $V^\lambda\leqslant 1$ for all $\lambda\in{\cal H}^*$.
\vskip 5pt
 In this paper, we will consider some properties of
${\cal L}$, basing on representations of the above type over the
Virasoro algebra.

The paper is arranged as follows. In Section 2, we first consider
all possible super-extensions of the Heisenberg-Virasoro type Lie
algebra. Our main result in this section is Theorem 2.1. As a
conclusion, we obtain that the Ramond $N=2$ superconformal algebra
is a special case of such super-extension. Then we study the modules
of the intermediate series  in the last section.

  Our main result is \\[5pt] {\bf Theorem 1.2.}\quad Any
indecomposable module of
  the intermediate series $V$ over the Ramond $N=2$ superconformal
  algebra is one of modules $RA_{a,b}$, $RA_\alpha$, $RA^\beta$,
  $RB_{a,b}$, $RB_\alpha$, $RB^\beta$,
  or one of their quotients
  for $a, b, \alpha, \beta\in{\ma C}$, where $RA_{a,b}$ is defined in
  (\ref{a50}), $RB_{a,b}$ is defined in (\ref{a71}), $RA_\alpha$ is defined in (\ref{b3}), $RA^\beta$ is defined in (\ref{b5}),
  (\ref{b7}), $RB_\alpha$ is defined in (\ref{c1}), $RB^\beta$ is defined in (\ref{c3}).
\vskip 10pt \cl{\bf 2. The structure of the Ramond $N=2$
superconformal algebras}
\def\theequation{2.\arabic{equation}}
\setcounter{equation}{0}\vskip5pt

Let $\tilde{\cal L}_{\overline{0}}={\rm span}_{\ma C}\{ L_m, H_n\mid
m,n\in{\ma Z}\}$ be a Heisenberg-Virasoro type algebra (only with
$1$-dimensional center ${\ma C} H_0$) with the following Lie
brackets: \BE\label{a1} [ L_m, L_n ]=(m-n)L_{m+n}, ~~ [ H_m, H_n]=0,
~~ [ L_m, H_n ]= -nH_{m+n}.
 \EE
Let us consider all possible super-extensions of the Lie algebra
$\tilde{\cal L}_{\overline{0}}$. First assume that $\tilde{\cal
L}_{\overline{1}}$ is an $\tilde{\cal L}_{\overline{0}}$-module of
intermediate series with basis $\{G_i\mid i\in{\ma
 Z}\}$ such that $\tilde{\cal
L}_{\overline{0}}\oplus \tilde{\cal L}_{\overline{1}}$ is a Lie
superalgebra. Then following \cite[Theorem 3.2]{LZ}, we can suppose
$$ [ L_m, G_i ]=(a-i+mb)G_{m+i}, ~~~ [ H_n, G_i ]=fG_{m+i},$$ for
some $a, b, f\in{\ma C}, f\not=0 $. Set $[ G_i, G_j
]=a_{ij}L_{i+j}+b_{ij}H_{i+j},$ for some $a_{ij},b_{ij}\in{\ma C}$.
Then from the equation \BC [ H_k, [ G_i, G_j ]]=[ [ H_k, G_i ], G_j
]+[ G_i, [ H_k, G_j ] ],\EC by letting $k=0$, we deduce $ 2f[ G_i,
G_j ]=0. $ Therefore, $[ G_i, G_j ]=0$ for all $i, j\in{\ma Z}.$
This is a trivial extension and not the thing we are interested in.
Hence we suppose that $\tilde{\cal L}_{\overline{1}}$ is a direct
sum of two $\tilde{\cal L}_{\overline{0}}$-modules of intermediate
series, with basis $\{G_i^{\pm}\mid i\in{\ma Z}\}$. Then we have the
following equations:
$$ [ L_m, G_i^+ ]=(a^+-i+mb^+)G_{m+i}^+, ~~~ [ H_n, G_i^+ ]=f_1G_{m+i}^+,$$
$$ [ L_m, G_i^- ]=(a^--i+mb^-)G_{m+i}^-, ~~~ [ H_n, G_i^- ]=f_2G_{m+i}^-,$$
where $a^+, a^-, b^+, b^-, f_1, f_2\in{\ma C}$, and $(f_1,
f_2)\ne(0,0)$. In order to get a nontrivial super-extension, we must
have $a^+=a^-$, denoted by $a$. Set
$$[ G_i^-, G_j^+ ]=a_{ij}L_{i+j}+b_{ij}H_{i+j}.$$ Since
$$ 0=[ H_0, [ G_i^-, G_j^+ ]]=(f_1+f_2)[ G_i^-, G_j^+ ],$$ we can
suppose that $f_1=1, f_2=-1$ (replacing $H_m$ by $f_1^{-1}H_m$).
Following \BE\label{a3} [ L_k, [ G_i^-, G_j^+ ]]=[ [ L_k, G_i^- ],
G_j^+ ]+[ G_i^-, [ L_k, G_j^+ ]], \EE and setting $k=0$ in
(\ref{a3}), we can get that
 $ a=0.$ Following (\ref{a3}), we have
\BE\label{a8}
(k-i-j)a_{i,j}&\!\!\!=\!\!\!&(-i+kb^-)a_{i+k,j}+(-j+kb^+)a_{i,k+j},\\
 (-i-j)b_{i,j}&\!\!\!=\!\!\!&(-i+kb^-)b_{i+k,j}+(-j+kb^+)b_{i,k+j}.\EE
 By \BE\label{a4}[ G_k^-, [ G_i^-, G_j^+ ]]+[ G_i^-, [
G_k^-, G_j^+ ]]=0, \EE we obtain $
b_{ij}+b_{kj}=a_{ij}(-k+(i+j)b^-)+a_{kj}(-i+(k+j)b^-).$ Setting
$k=i$ gives \BE\label{a5} b_{ij}=a_{ij}(-i+(i+j)b^-).\EE Replacing
$G_k^-, G_i^-, G_j^+ $ respectively by $G_j^+, G_j^+, G_i^-$ in
(\ref{a4}) gives
 $$
 b_{ij}=-a_{ij}(-j+(i+j)b^+). $$ Comparing it with (\ref{a5}), we can get that $$b^++b^-=1.$$ Following
$ [ H_k, [ G_i^-, G_j^+ ]]=[ [ H_k, G_i^- ], G_j^+ ]+[ G_i^-, [ H_k,
G_j^+ ] ],$ we have that \BC a_{k+i,j}=a_{i,k+j}, \ \ \ ka_{i,j}=
b_{i,k+j}-b_{k+i,j} ~~~\mbox{for all} ~~ i, j, k\in{\ma Z}.\EC
Letting $i=0, k=i$ gives \BE\label{a7} a_{i,j}=a_{0,i+j} ~~\mbox{for
all} ~ i, j\in{\ma Z}.\EE   Taking $i=0$ in (\ref{a8}), and by
(\ref{a7}), we obtain that  \BE\label{a9}(k-j)a_{0,j}=kb^-
a_{k,j}+(-j+kb^+)a_{0,k+j}=(k-j)a_{0,k+j}.\EE
 Let
 $a_{0,0}=d\in{\ma C}$, then by (\ref{a7}) and (\ref{a9}), we have that $a_{ij}=d~ \mbox{for all } i, j\in{\ma Z}.$
Since $b^++b^-=1$, we set $b^-=b$, then $$b^+=1-b, ~~\mbox{and} ~~
b_{ij}=d(-i+(i+j)b).$$ We can set $d=1$ (replace $G^\pm_i$ by
$\frac{1}{\sqrt{d}}G^\pm_i$), then $\tilde{\cal L}=\tilde{\cal
L}_{\overline{0}}\oplus\tilde{\cal L}_{\overline{1}}$ is a
superalgebra with (\ref{a1}) and the following Lie brackets:
\BE\label{a21}\begin{array}{lllllll} [ L_m, G_n^+ ]=
(-n+m(1-b))G_{m+n}^+, \ \ \ \ [ L_m, G_n^- ]= (-n+mb)G_{m+n}^- ,
\\ [7pt]
 [ H_m, G_n^{\pm} ] ={\pm}G_{m+n}^{\pm}, \ \ \ \ \ \ \
 \ ~~~~~~~~~~~~~~~
 [ G_m^\pm, G_n^\pm ]= 0, \\[7pt]
 [ G_m^-, G_n^+ ]= L_{m+n}+(-m+(m+n)b)H_{m+n},
 \end{array}\EE where $m, n\in{\ma Z}.$ Obviously, $\tilde{{\cal L}}$ is ${\mathbb Z}$-graded:
\BC
 \tilde{{\cal L}}= \OPLUS{n\in{\mathbb Z}}\tilde{{\cal L}}_n, \ \ \
 \
\tilde{{\cal L}}_n =\{x\in \tilde{\cal L}\mid [L_0,
x]=nx\}=\mbox{span}_{\mathbb C} \{L_{-n}, H_{-n}, G^{\pm}_{-n}
 \}.\EC
\indent
  Now we consider the central extension of
$\tilde{\cal L}$. Suppose $\varphi: \tilde{\cal L}\times \tilde{\cal
L}\rightarrow {\ma C}$ is a 2-cocycle of $\tilde{\cal L}$, we define
a linear map $f: \tilde{\cal L}\rightarrow{\ma C}$ as follows:
\BE\label{a12}\begin{array}{lll} f({\tilde{\cal
L}_i})={\frac{1}{i}}\varphi(L_0, \tilde{\cal L}_i),
~~ i\not=0, & f(L_0)={\frac{1}{2}}\varphi(L_1, L_{-1}), \\[0.4 cm]
f(H_0)=\varphi(L_1, H_{-1}),& f(G^\pm_{0})=
\pm\varphi(H_1,G_{-1}^\pm).
\end{array}\EE If we define another  2-cocycle of $\tilde{\cal L}$,
$\psi: \tilde{\cal L}\times \tilde{\cal L}\rightarrow {\ma C}$
satisfying $\psi=\varphi-\varphi_f$, where
$\varphi_f(x,y)=f([x,y])$, then \BC \psi(L_0, \tilde{\cal
L}_{i})&\!\!\!=\!\!\!& \varphi(L_0, \tilde{\cal
L}_{i})-\varphi_f(L_0, \tilde{\cal L}_{i}) =  0, ~~~ i\not=0, \\
[0.2 cm] 2\psi(L_0,L_0)&\!\!\!=\!\!\!& \psi([L_1, L_{-1}], L_0) =
0,\EC Similarly, we have that $ \psi(L_0, H_0)= \psi(L_0,
G_{0}^\pm)= 0.$ Thus \BE\label{a10}\psi(L_0, \tilde{\cal L}_i)=0 ~~
\mbox{ for all }
 i\in{\ma Z}.\EE Furthermore, \BC  i\psi(\tilde{\cal L}_i,
\tilde{\cal L}_j) &\!\!\!=\!\!\!&
\psi([L_0, \tilde{\cal L}_i], \tilde{\cal L}_j)
=
\psi(L_0, [\tilde{\cal L}_i, \tilde{\cal
L}_j])-\psi(\tilde{\cal
L}_i, [L_0, \tilde{\cal L}_j])
=
-j\psi(\tilde{\cal L}_i, \tilde{\cal L}_j), \EC
therefore, \BE\label{a11} \psi(\tilde{\cal L}_i, \tilde{\cal L}_j)=0
~~~ \mbox{if} ~~ i+j\not=0 .\EE Now let us consider $\psi(L_i,
L_{-i})$. It follows (\ref{a12}) that \BC \psi(L_1, L_{-1})=
\varphi(L_1, L_{-1})-\varphi_f(L_1, L_{-1})= 0. \EC Then \BC
(i-2)\psi(L_i, L_{-i}) = \psi([L_{i-1}, L_{1}], L_{-i})
  = (i+1)\psi(L_{i-1}, L_{-i+1}). \EC
Set $\psi(L_i, L_{-i})=l_i$, we obtain $l_i=\frac{i+1}{i-2}l_{i-1}$
for $i\not= 2,$  i.e.,
 $$\mbox{$l_i=\frac{i^3-i}{6}l_2 ~~~\mbox{for all } ~~ i\geqslant 3.$}$$ Then we can rewrite  $\psi(L_i, L_{-i})$ as follows
 \BE\label{a13} \psi(L_i, L_{-i})=\mbox{$\frac{i^3-i}{6}$}c_L ~~~ \mbox{for all } i\in{\ma Z},\EE where $c_L\in{\ma
 C}$.  Similar to the
argument about $\psi(L_i, L_{-i})$, we can obtain that
 \BE
&\!\!\!\!\!\!\!\!\!\!\!\!\!\!\!\!& \psi(H_i, H_{-i})=ic_H, \ \ \ \ \
 \label{a15}\psi(L_i,
 H_{-i})=\mbox{$\frac{i(i-1)}{2}$}c_{HL},\\[5pt]
&\!\!\!\!\!\!\!\!\!\!\!\!\!\!\!\!& \label{a16} \psi(H_i,
G_{-i}^\pm)=\psi(L_i, G_{-i}^\pm)=\psi(G_i^\pm, G_{-i}^\pm)=0
~~\mbox{ for all } i\in{\ma Z}, \EE where $c_H, c_{HL}\in{\ma C}$.
  Finally, we consider $\psi(G_i^-, G_{-i}^+)$. Set
$\psi(G_0^-, G_0^+)=c_G$, by (\ref{a21}) 
and
(\ref{a15}), we have \BC
\psi(G_1^-, G_{-1}^+)&\!\!\!=\!\!\!& \psi([G_0^-, H_1], G_{-1}^+)\\
&\!\!\!=\!\!\!& \psi(G_0^-,
[H_1, G_{-1}^+])-\psi(H_1, [G_0^-, G_{-1}^+])
=
c_G+c_{HL}+bc_H, \EC and $ \psi(G_1^-, G_{-1}^+) = \psi(G_1^-,
[H_{-1}, G_0^+]) = c_G+(1-b)c_H. $ Hence
$c_G+c_{HL}+bc_H=c_G+(1-b)c_H$,
i.e.,\BE\label{a19}c_{HL}=(1-2b)c_H.\EE Then
\BE\label{a20}\begin{array}{lll} \psi(G_i^-,
G_{-i}^+)&\!\!\!=\!\!\!& \psi([G_0^-,
H_i], G_{-i}^+)\\[0.2 cm]
&\!\!\!=\!\!\!&\psi(G_0^-, [H_i, G_{-i}^+])-\psi(H_i, [G_0^-, G_{-i}^+])
=
c_G+\big({
\frac{i(i+1)}{2}}-ib\big)c_H.
\end{array}\EE Note that  \BC [L_i, [G_j^-, G_k^+]]&\!\!\!=\!\!\!&
[[L_i, G_j^-], G_k^+]+[G_j^-, [L_i, G_k^+]]. \EC  If
 we suppose  $i+j+k=0$, then  \BE\label{a22}\begin{array}{lll} &\!\!\!\!\!\!\!\!\!\!\!\! &
\frac{i^3-i}{6}c_L+(-j+(k+j)b)\frac{i(i-1)}{2}c_{HL} \\[4pt]
  &\!\!\!\!\!\!\!\!\!\!\!\!&
  =(-j+ib)\big(c_G+(\frac{-k(1-k)}{2}+kb)c_H\big)
  +(-k+i(1-b))
  \big(c_G+(\frac{j(j+1)}{2}-jb)c_H\big).
 \end{array}\EE
By (\ref{a19}), and setting $j=0$ in (\ref{a22}), we have that
\BE\label{a25} \mbox{$\frac{i^2-1}{6}$}c_L-2c_G=i^2(b-b^2)c_H. \EE
Letting $i=1$ in (\ref{a25}), we can obtain that  \BE\label{c_G}
c_G=\mbox{$\frac{b^2-b}{2}$}c_H, ~~~~ c_L=6(b-b^2)c_H.\EE Then we
have the following theorem: \\ [5pt] {\bf Theorem 2.1.}\quad The
possible nontrivial super-extensions of the Heisenberg-Virasoro type
algebra (\ref{a1}) are the following superalgebras:
$$ \hat{\cal L}=\mbox{span}_{\ma C}\{L_i, H_j, G^\pm_k, c_H \mid c_H\in{\ma C}, i, j, k\in{\ma Z} \},
$$ where $c_H$ is a central element and the following relations hold:
\BE\label{a}\begin{array}{lll}
 [ L_i, L_j ]=
 (i-j)L_{i+j}+(i^3-i)(b-b^2)c_H\delta_{i+j,0}, \\ [5pt]
   [ L_i, H_j ] =-jH_{i+j}+\frac{i(i-1)}{2}(1-2b)c_{H}\delta_{i+j,0}, \ \
   [ H_i, H_j ]= ic_H \delta_{i+j,0}, \\ [5pt]
    [ L_i, G_j^+ ]= (-j+i(1-b))G_{i+j}^+, \  \ \ \ \ \ \ \ \ \ \ \ \
    \ \ \
[ L_i, G_j^- ] = (-j+ib)G_{i+j}^- , \\[5pt]
   [ H_i, G_j^{\pm} ]=\pm G_{i+j}^{\pm}, \ \ \ \ \ \ \ \ \ \ \ \ \ \
   \ \ \ \ \ \ \ \ \ \ \ \ \ \ \ \ \ \ \
    [ G_i^+, G_j^+ ]= [ G_i^-, G_j^- ]= 0, \\ [5pt]
    [ G_i^-, G_j^+ ]=
  L_{i+j}+(-i+(i+j)b)H_{i+j}+\frac{i(i+1-2b)+b^2-b}{2}c_H\delta_{i+j,0}.\end{array} \EE
   \vskip 7pt
If $b=\frac{1}{2}$, then $\hat{\cal L}={\cal L}$. That is to say,
Ramond $N=2$ superconformal algebra ${\cal L}$ is a special case of
$\hat{\cal L}$.
\vskip 9pt

\cl{\bf 3. The modules of intermediate series over ${\cal L}$}
\def\theequation{3.\arabic{equation}}
\setcounter{equation}{0} \vskip4pt
 \noindent{\bf \S3.0}\quad
 Let $V=V_0\oplus V_1$ be any indecomposable ${\cal L}$-module with
 dim$V^{\lambda}_{\alpha}\leqslant 1$ for all $\lambda\in{\cal H}^*, \alpha\in{\mathbb Z}/2{\mathbb
 Z}$, where $$V^{\lambda}_{\alpha}=\{v\in V_{\alpha}\mid L_0\cdot v=\lambda(L_0)v, H_0\cdot v=\lambda(H_0)v
 \}.$$ We also have the following:
\BE\label{a30} V=\Big(\OPLUS{k\in{\mathbb
Z}}V_0^{a+k}\Big)\oplus\Big(\OPLUS{k\in{\mathbb Z}}V_1^{a+k}\Big).
\EE
 One sees that $c$ acts trivially on $V$ (see, e.g., \cite{IL, CA}). So
 we can omit $c$ in (\ref{a23}).

 Now we consider all possibilities of $V_0$ and $V_1$ case by case
 below. Let us recall the definition of $Vir$-modules $A_{a,b}, A(\alpha),
 B(\beta)$ (see \cite{IL}). They all have a basis $\{x_i\mid i\in{\mathbb
 Z}\}$ such that for $i, j\in{\mathbb Z}$,
 \BE\label{a31}\begin{array}{llll}
  A_{a,b}:&  L_ix_j=(a-j+ib)x_{i+j}. & &   \\ [0.2 cm]
  A(\alpha):&  L_ix_j=-(i+j)x_{i+j}, ~ j\not=0, &
L_ix_0 = -i(1+(i+1)\alpha)x_i. \\ [0.2 cm]
  B(\beta):&  L_ix_j=-jx_{i+j}, ~ i+j\not=0, & L_ix_{-i} = i(1+(i+1)\beta)x_0.
 \end{array}\EE
 \\
 {\bf \S3.1}\quad Suppose both of $V_0, V_1$ have the form $A_{a,b}, a, b\in {\mathbb
 C}$. Then we choose a basis $\{x_i\mid i\in{\mathbb Z}\}$ of $V_0$
 and a basis $\{y_j\mid j\in{\mathbb Z}\}$ of $V_1$ such that
\BE\label{a32} L_ix_j=(a-j+ib)x_{i+j}, ~~~ L_iy_j=(a'-j+ib')y_{i+j},
\\ [0.1 cm] \label{a33} H_ix_j=f_{ij}x_{i+j}, \hspace{2cm}
H_iy_j=f'_{ij}y_{i+j},\hspace{1.8cm}
\\ [0.1 cm]
\label{a34} G^{\pm}_ix_j=a_{ij}^{\pm}y_{i+j}, \hspace{1.9cm}
G^{\pm}_iy_j=b_{ij}^{\pm}x_{i+j},\hspace{1.8cm} \EE where $a, \ a',\
b,\ b',\ f_{ij},\ f'_{ij},\ a^\pm_{ij},\ b^\pm_{ij}\in{\ma C}$.
 We have $a=a'$ by
applying $L_0$ to the first formula of (\ref{a34}). By (\ref{a23}),
we have \BE\label{a35}
(\mbox{$\frac{i}{2}$}-j)a_{i+j,k}^{\pm}=(a-(k+j)+ib')a_{j,k}^{\pm}-(a-k+ib)a_{j,i+k}^{\pm},
\EE and
 $ 2(a-k+(i+j)b)-(i-j)f_{i+j,k}=a_{j,k}^{+}b_{i,k+j}^{-}+a_{i,k}^-b_{j,k+i}^+.$
Let $i=j$, we get \BE\label{a36}
a_{i,k}^{+}b_{i,k+i}^{-}+a_{i,k}^-b_{i,k+i}^+=2(a-k+2ib). \EE From
(\ref{a36}), we know that for all fixed $i, k\in{\mathbb Z}$,
\BE\label{a37} a_{i,j}^+=a_{i,j}^-=0,~ b_{k,l}^+=b_{k,l}^-=0
~\mbox{only for finitely many $j$ and finitely many $l$}.\EE
Applying $ [ G_i^\pm, G_j^\pm ]=0$ to $x_k$ gives $ a_{i,j}^\pm
b_{i,k+j}^\pm+a_{i,k}^\pm b_{j,k+i}^\pm=0.$ Letting $i=j$ gives
\BE\label{a46} a_{i,k}^\pm b_{i,k+i}^\pm=0. \EE Therefore, by
(\ref{a36}) and (\ref{a46}), for any $i, k\in{\ma Z}$ with
$a-k+2ib\not=0$, we have \BE\label{a48} a_{ik}^+a_{ik}^-=0, ~~
(a_{ik}^+)^2+(a_{ik}^-)^2\not=0, ~~ (\mbox{similar relations for} ~~
b_{ik}^\pm).\EE
 For
convenience, we omit the superscript ``\,$^{\pm}$\,'' in
$a^\pm_{ij}$. Let $i=2j$ and $i=-2j$ in (\ref{a35}) respectively, we
get \BE\label{a38}
(a-k+2jb)a_{j,k+2j}&\!\!\!=\!\!\!& (a-(k+j)+2jb')a_{j,k}, \\[4pt]
\label{a39} -2ja_{-j,k}&\!\!\!=\!\!\!&
(a-(k+j)-2jb')a_{j,k}-(a-k-2jb)a_{j,k-2j}. \EE Multiplying
(\ref{a39}) by $a-(k-j)+2jb'$ and replacing the last term by
(\ref{a38}), we get
\begin{eqnarray*}
&\!\!\!\!\!\!\!\!\!\!\!\! & -2j(a-(k-j)+2jb')a_{-j,k}\\[4pt]
\!\!\!\!\!\!\!\!\!\!\!\!
&\!\!\!\!\!\!\!\!\!\!\!\!&=(a\!-\!(k\!-\!j)\!+\!2jb')(a\!-\!(k\!+\!j)\!-\!2jb')a_{j,k}
\! -\!(a\!-\!k\!-\!2jb)(a\!-\!(k\!-\!2j)\!+\!2jb)a_{j,k}\\[4pt]
\!\!\!\!\!\!\!\!\!\!\!\! &\!\!\!\!\!\!\!\!\!\!\!\!&=
2j(a-(k-j)+2jb'+2jt)a_{j,k},
 \end{eqnarray*}
where $t=b'^2-(b+\frac{1}{2})^2$. Similarly, let $j=-j$, $i=2j$ and
$j=-j$, $i=-2j$ in (\ref{a35}), we can obtain that
\begin{eqnarray*}
2j(a-(k+j)-2jb')a_{j,k}=-2j(a-(k+j)-2jb'-2jt)a_{-j,k}.
 \end{eqnarray*}
It follows that
$$\begin{array}{ll}\big( (a-(k+j)-2jb')(a-(k-j)+2jb')\\[3pt]
-((a-k-j-2jb')-2jt)((a-k+j+2jb')+2jt)\big)a_{jk}  =0,\end{array}$$
which gives \BE\label{a40} 4j^2t(t+2b'+1)a_{jk}=0. \EE
 By (\ref{a37}), there at least exists
one $k_0$ such that $a_{1,k_0}^+\not=0$, or $a_{1,k_0}^-\not=0$. (If
$a^\pm_{1,k}=0$ for all $k$, by letting $j=1$ in (\ref{a35}) we get
$a_{ik}^\pm=0$
 for $i, k\in{\mathbb Z}, i\not=2$. By letting $i=j=1$ in (\ref{a35}), we get
 $a^\pm_{2,k}=0.$)
Thus it follows from (\ref{a40}) that $$\mbox{$
b'=\pm(b+\frac{1}{2}), ~~\mbox{or}~~ b'=-1\pm(b+\frac{1}{2}). $}$$
\\[5pt]{\bf Case 1.}\quad $b'=b+\frac{1}{2}$.\vskip 4pt

First suppose $a-k+2b\not=0$ for all $k\in{\mathbb Z}$. Letting
$j=1$ in (\ref{a38}), we obtain (again we omit the superscript
``\,$^\pm$\,'' in $x_0,\,x_1$ for the time being)
\begin{eqnarray*} a_{1,k}=\left \{
\begin{array}{ll}
x_0, & k\mbox{ is even}, \\
x_1, & k \mbox{ is odd}.  \end{array}\right.
\end{eqnarray*}
 Let
$j=1$, and let $i, k$ be odd in (\ref{a35}), then
\begin{eqnarray*}
(\mbox{$\frac{i}{2}$}-1)a_{i+1,k}&\!\!\!=\!\!\!&(a-(k+1)+ib')a_{1,k}-(a-k+ib)a_{1,k+i} \\
&\!\!\!=\!\!\!& (a-(k+1)+ib')x_1-(a-k+ib)x_0.
\end{eqnarray*}
By (\ref{a38}) and $b'=b+\frac{1}{2}$, we also have
\begin{eqnarray*}
(\mbox{$\frac{i}{2}$}-1)a_{i+1,k}&\!\!\!=\!\!\!& (\mbox{$\frac{i}{2}$}-1)a_{i+1,k+2(i+1)}\\
&\!\!\!=\!\!\!& (a-(k+2(i+1)+1)+ib')x_1-(a-(k+2(i+1))+ib)x_0.
\end{eqnarray*}
Obviously, we get $x_0=x_1$. Similar to the arguments after
(\ref{a40}), we have
\begin{equation}\label{a-i-j=d-i-j}a_{ij}^\pm =d_1^{\pm} ~~~
\mbox{for all} ~~ i,j\in{\ma Z},
\end{equation} where $d_1^\pm$ are constants, and by (\ref{a48}),
\begin{equation}\label{d-i-jnot=}
d_1^+d_1^-=0, ~~~ (d_1^+)^2+(d_1^-)^2\not=0. \end{equation} Now we
suppose that $a-k'+2b=0$ for some $k'\in{\mathbb Z}$. It follows
from (\ref{a38}) that $ (a-k+2b)a_{1,k}=(a-k+2b)a_{1,k+2}. $ Then
 $$ \begin{array}{ll} a_{1k}=\left \{
\begin{array}{ll}
x_0, & k>k', ~ k \mbox{ is  even}, \\
x_1, & k>k', ~ k  \mbox{ is  odd}.  \end{array}\right. &
a_{1k}=\left \{
\begin{array}{ll}
y_0, & k\leqslant k', ~  k \mbox{ is  even}, \\
y_1, & k\leqslant k', ~ k \mbox{ is  odd}.
\end{array}\right.\end{array}
$$
By (\ref{a35}), we get
\begin{eqnarray*} a_{ik}=\left \{
\begin{array}{ll}
x_0, & k>k', ~ k+i-1> k', ~\mbox{and} ~ k , ~i-1 \mbox{  are even}, \\
x_1, & k>k', ~ k+i-1> k', ~\mbox{and} ~ k , ~i-1 \mbox{  are  odd},
\end{array}\right.
\end{eqnarray*}
and
\begin{eqnarray*} a_{ik}=\left \{
\begin{array}{ll}
y_0, & k\leqslant k', ~ k+i-1\leqslant k', ~\mbox{and} ~ k , ~i-1 \mbox{ are  even}, \\
y_1, & k\leqslant k', ~ k+i-1\leqslant k', ~\mbox{and} ~ k , ~i-1
\mbox{ are odd}.
\end{array}\right.
\end{eqnarray*}
Now choose some $k, j\in{\mathbb Z}$,  such that $a-k+2jb\not=0, ~
k\leqslant k', ~ k+j-1\leqslant k'$, $k+2j>k', ~ k+3j-1> k'$, and
one of $k$ and $j$ is even, and another is odd. Then by (\ref{a38}),
we have $$ (a-k+2jb)a_{j,k}=(a-k+2jb)a_{j,k+2j}. $$ Therefore, $$
x_0=y_0, ~~~ x_1=y_1.$$ Similar to the argument above, we again have
(\ref{a-i-j=d-i-j}) and (\ref{d-i-jnot=}).
\\[4pt]
{\bf Case 2.}\quad $b'=-(b+\frac{1}{2})$. \vskip4pt

Let $i=2j$ and $i=-2j$ respectively in (\ref{a35}), we have
\BE\label{a41}\!\!\!\!&\!\!\!\!&
(a-k-2j-2jb)a_{j,k}=(a-k+2jb)a_{j,k+2j},\\\label{a42}
\!\!\!\!&\!\!\!\!&
-2ja_{-j,k}=(a-k+2jb)a_{j,k}-(a-k-2jb)a_{j,k-2j}=-2ja_{j,k}, \EE
where the second equality of (\ref{a42}) follows from (\ref{a41}) by
replaced $k$ by $k-2j$. Hence, $a_{j,k}=a_{-j,k}$ for all $k,
j\in{\mathbb Z}$. Using it, again by (\ref{a35}), we deduce that
\begin{eqnarray*}
 & & (a-(k+j)+ib')a_{j,k}-(a-k+ib)a_{j,k+i}
 =
 (\mbox{$\frac{i}{2}$}-j)a_{i+j,k}\\
&& =(\mbox{$\frac{i}{2}$}-j)a_{-i-j,k}
 = -(a-(k-j)-ib')a_{-j,k}-(a-k-ib)a_{-j,k-i}.
 \end{eqnarray*}
Hence \BE\label{a43}
(a-k-ib)a_{j,k-i}-2(a-k)a_{j,k}+(a-k+ib)a_{j,k+i}=0. \EE Let $j=1$
in (\ref{a41}) and then replace $k$ by $k+2$ in the new equality, we
can obtain
$$\mbox{$ a_{1,k+4}=\frac{(a-k-4-2b)(a-k-2-2b)}{(a-k-2+2b)(a-k+2b)}a_{1,k}.$} $$
Similarly, we can have a formula for $a_{1,k-4}$. Then let $i=4$ in
(\ref{a43}), we have
\begin{eqnarray*}& &
\mbox{$\big((a-k-4b)\frac{(a-k+4+2b)(a-k+2+2b)}{(a-k+2-2b)(a-k-2b)}-2(a-k)$}
\\[4pt]
 & & \mbox{$+(a-k+4b)
 \frac{(a-k-4-2b)(a-k-2-2b)}{(a-k-2+2b)(a-k+2b)}\big)a_{1,k}=0.$}
\end{eqnarray*}
By (\ref{a38}) and the discussion after (\ref{a40}), we know  the
coefficient of $a_{1,k}$ must be zero. We obtain  $b=-1$ or
$-\frac{1}{2}$. Note that the case of $b=-\frac{1}{2}, b'=0$ is
contained in Case 1. So we only need to consider the case of $b=-1,
b'=\frac{1}{2}$.

Let $j=1$ in (\ref{a38}), then  $(a-k)a_{1,k}$ is a constant for all
even $k$ or all odd $k$. We suppose that
\begin{eqnarray*} (a-k)a_{1,k}=\left \{
\begin{array}{ll}
x_0, & k\mbox{ is even}, \\
x_1, & k \mbox{ is odd}.  \end{array}\right.
\end{eqnarray*}
Let $i=j=1$ in (\ref{a43}), we can obtain that $x_0=x_1$. That is to
say, $(a-k)a_{1,k}$ is a constant for all $k\in{\ma Z}$.
 If $a-k_1=0$ for some $k_1\in{\mathbb Z}$,
then $ a_{1,k}=\frac{a-k_1}{a-k}a_{1,k_1}=0 ~\mbox{for all} ~
k\not=k_1, $
a contradiction with (\ref{a37}). Thus $a-k\not=0$ for all
$k\in{\mathbb Z}$, i.e., $a\not\in{\ma Z}$. Denote
$(a-k)a_{1,k}^\pm$ by $d_2^\pm$, where $d_2^\pm\in{\mathbb C}$ are
constants. By (\ref{a35}), let $j=1$, we get
$$ (\mbox{$\frac{i}{2}$}-1)a_{i+1,k}=(a-(k+1)+\mbox{$\frac{i}{2}$})a_{1,k}-(a-k-1)a_{1,k+1},
$$ then $ a_{i,k}=a_{1,k} ~\mbox{for all} ~ i\in{\mathbb Z}, i\not=3. $
Again by (\ref{a35}), let $j=2$,  $i=1$, we have  $$
-\mbox{$\frac{3}{2}$}a_{3,k}=(a-(k+2)+\mbox{$\frac{1}{2}$})a_{2,k}-(a-k-1)a_{2,k+1},
$$ then $a_{3,k}=a_{1,k}.$
Therefore, $ a_{ij}^\pm=(a-j)^{-1}d_2^\pm ~\mbox{for all} ~ i,
j\in{\mathbb Z}, $ and by (\ref{a48}), $$ d_2^+d_2^-=0, ~~~
(d_2^+)^2+(d_2^-)^2\not=0.
$$
{\bf Case 3.}\quad $b'=-b-\frac{3}{2}$. \vskip4pt  Following the
arguments in Case 2, we have $ 2ja_{j,k}=2ja_{-j, k-2j}, $ and
\BE\label{b4} \!\!\!\!\!\!\!\!\!\!&\!\!\!\!\!\!\!\!\!\!&
(a-k-j+ib')a_{j,k}+(a-k-j-2i-ib')a_{j,k+2i}=2(a-k-i-j)a_{j,k+i}, \\
[0.2 cm]
 \!\!\!\!\!\!\!\!\!\!&\!\!\!\!\!\!\!\!\!\!&
\begin{array}{lll}\mbox{$
\big((a-k-1+4b')\frac{(a-k-3-2b')(a-k-5-2b')}{(a-k-1+2b')(a-k-3+2b')}-2(a-k-5)
$} \\ [0.2 cm]
 \mbox{$\ \ \ \ \ \ \ \ \ \ \ \ \ \ \ \ \ \ \ \ \
 +(a-k-9-4b')\frac{(a-k-7+2b')(a-k-5+2b')}{(a-k-9-2b')(a-k-7-2b')}\big)a_{1,k+4}=0.$}
 \end{array}\EE
 Then we obtain that $b=-\frac{3}{2}$ or $-1$. The case of $b=-1\
 (b'=-\frac{1}{2})$ is contained in Case 1. So we consider the
 case of $b=-\frac{3}{2}, b'=0$. By (\ref{a38}), let $j=1$, we have
 $$ (a-k-1)a_{1,k}=(a-k-3)a_{1,k+2}, $$ then
 $(a-k-1)a_{1,k}$ is a constant for all
even $k$ or all odd $k$. By (\ref{b4}), similar to the argument in
Case 2, we can obtain that $(a-k-1)a_{1,k}$ is a constant for all
$k\in{\ma Z}$. If $a-k_1-1=0$ for some $k_1$, then
$a_{1,k}=\frac{a-k_1-1}{a-k-1}a_{1,k_1}=0$ for all $k\ne k_1$. Also
a contradiction with (\ref{a37}). Therefore, $a\not\in{\ma Z}$. Now
we
 denote $(a-k-1)a_{1,k}$ by $d_3$, i.e. $a_{1,k}=(a-k-1)^{-1}d_3$,
for all $k\in{\ma Z}$. Let $j=1$ in (\ref{a35}), then we have that
$$\mbox{$ (\frac{i}{2}-1)a_{i+1,k}=(a-k-1)a_{1,k}-(a-k-\frac{3}{2}i)a_{1,k+i}.$} $$
Then $$ a_{i,k}=a_{1,k+i-1}=(a-k-i)^{-1}d_3 ~~ \mbox{for } i\ne 3.$$
Let $i=1, j=2$ in (\ref{a35}), we get $-\frac{3}{2}a_{3,k}=
(a-k-2)a_{2,k}-(a-k-\frac{3}{2})a_{2,k+1}=-\frac{3}{2}a_{2,k+1}$,
then $a_{3,k}=(a-k-3)^{-1}d_3$.
 Therefore, $ a_{ij}^\pm=(a-i-j)^{-1}d_3^\pm ~\mbox{for all} ~ i,
j\in{\mathbb Z}, $ and by (\ref{a48}), $$ d_3^+d_3^-=0, ~~~
(d_3^+)^2+(d_3^-)^2\not=0.
$$
\\{\bf Case 4.}\quad $b'=b-\frac{1}{2}$.\vskip 4pt
Note that if we act $(\frac{i}{2}-j)G_{i+j}^\pm=[ L_i, G_j^\pm ]$ on
$y_k$, we can obtain that \BE\label{a44}
(\mbox{$\frac{i}{2}$}-j)b_{i+j,k}^\pm=(a-(k+j)+ib)b_{j,k}^\pm-(a-k+ib')b_{j,i+k}^\pm.
\EE Similar to the discussion in case 1 $(b=b'+\frac{1}{2})$, we
have $ b_{ij}^\pm=d_1'^\pm,\,i, j\in{\ma Z}$ for some
$d_1'^\pm\in{\ma C}$, and by (\ref{a48}),
$$ d_1'^+d_1'^-=0, ~~~ (d_1'^+)^2+(d_1'^-)^2\not=0.
$$
Then it follows from (\ref{a36}) that
$$ a_{ij}^\pm=2(a-k+2ib)(d'^\mp_1)^{-1}. $$
Until now, we get that \begin{eqnarray*} a_{i j}^\pm=\left\{
\begin{array}{ll}
d_1^\pm, & b'=b+\frac{1}{2}, \\[4pt]
(a-j)^{-1}d_2^\pm, & b'=\frac{1}{2}, ~~ b=-1 ,\\[4pt]
(a-i-j)^{-1}d_3^\pm, & b'=0, ~~ b=-\frac{3}{2},\\[4pt]
2(a-j+2ib)(d'^\mp_1)^{-1}, & b'=b-\frac{1}{2}.
\end{array}\right.\end{eqnarray*}
Again using (\ref{a42}), following the same arguments about
$a_{i,j}^\pm$, we have that
\begin{eqnarray*} b_{ij}^\pm=\left\{
\begin{array}{ll}
d'^\pm_1, & b=b'+\frac{1}{2}, \\[4pt]
(a-j)^{-1}d'^\pm_2, & b=\frac{1}{2}, ~~ b'=-1 ,\\[4pt]
(a-i-j)^{-1}d'^\pm_3, & b=0, ~~ b'=-\frac{3}{2},\\[4pt]
2(a-j+2i(b+\frac{1}{2}))(d^\mp_1)^{-1}, & b=b'-\frac{1}{2}.
\end{array}\right.\end{eqnarray*}
Obviously only the following two cases can occur: \BE\label{a45}
\begin{array}{lll}
a_{ij}^\pm= d^\pm, & b_{ij}^\pm =
 2(a-j+2i(b+\frac{1}{2}))(d^\mp)^{-1}, &
 b'=b+\frac{1}{2}. \\[4pt]
 b_{ij}^\pm =d'^\pm, & a_{ij}^\pm =
 2(a-j+2ib)(d'^\mp)^{-1}, &  b=b'+\frac{1}{2},
\end{array}\EE
for some $d^\pm,d'^\pm\in{\ma C}$. Together with (\ref{a48}), we
obtain that
$$\begin{array}{l} (d^+)^2+(d^-)^2=0, ~~~ d^+ d^-=0,\mbox{ \ and }
\\[4pt] (d'^+)^2+(d'^-)^2=0, ~~~ d'^+ d'^-=0.
\end{array}$$ By rescaling basis $\{y_i\mid i\in{\ma Z}\}$ (or $\{x_i\mid
i\in{\ma Z}\}$) if necessary, we can suppose $d^\pm=1$ (or $
d'^\pm=1$). Then we rewrite (\ref{a45}):
\BE\label{a47}\begin{array}{llllllll} a_{ij}^+\ (\mbox{or} ~
a_{ij}^-) &\!\!\!=\!\!\!& 1, & b_{ij}^-\ (\mbox{or} ~ b_{ij}^+)
&\!\!\!=\!\!\!&
 2(a-j+2i(b+\frac{1}{2})), &
 b'=b+\frac{1}{2}. \\[4pt]
 b_{ij}^+\ (\mbox{or} ~ b_{ij}^-) &\!\!\!=\!\!\!& 1, & a_{ij}^-\ (\mbox{or} ~ a_{ij}^+) &
 = &
 2(a-j+2ib), &  b'=b-\frac{1}{2}.\end{array}
\EE

      Now we consider one of the cases of (\ref{a47}):
$$a_{ij}^+= 1, \ b_{ij}^- =
 2(a-j+2i(b+\mbox{$\frac{1}{2}$})),\ a_{ij}^-=b_{ij}^+=0,\ b'=b+\mbox{$\frac{1}{2}$}.$$  We want to determine the action of
$H_i$ on $V$. Set $f_{0 0}=f, f'_{0 0}=f'$. Since $[H_0,
G_i^+]=G_i^+$, we have that $f'=f+1$. Similar to the arguments of
\cite{LZ}, we get the following cases.
\\ [0.2 cm]
{\bf Case 5.}\quad $a-b, a-b'\not\in{\ma Z}$. We have some cases as
follows:
 \BC\begin{array}{llll} \mbox{Subcase
5.1.} &  f_{ij}=f, &  f_{ij}'=f+1; &  b'=b+\frac{1}{2}.\\
[0.3 cm] \mbox{Subcase 5.2.} &
f_{ij}=\frac{a-j}{a-i-j}f, & f_{ij}'=f+1; & b=0, ~ b'=\frac{1}{2}.\\
[0.3 cm] \mbox{Subcase 5.3.}& f_{ij}=\frac{a-i-j}{a-j}f, &
f_{ij}'=f+1; & b=-1, ~ b'=-\frac{1}{2}.\\ [0.3 cm] \mbox{Subcase
5.4.}&  f_{ij}=f, & f_{ij}'=\frac{a-j}{a-i-j}(f+1); &
b=-\frac{1}{2}, ~ b'=0.\\ [0.3 cm]\mbox{Subcase 5.5.}& f_{ij}=f, &
f_{ij}'=\frac{a-i-j}{a-j}(f+1); & b=-\frac{3}{2}, ~ b'=-1.
\end{array}\EC   Note that \BE\label{a49}[H_i,
G_j^+]=G_{i+j}^+. \EE Acting it on $x_k$, we compare the
coefficients on the two sides, then we get contradictions for
Subcases 5.2--5.5. And in Subcase 5.1, follows $[ H_i, G_j^-]\cdot
y_k=-G_{i+j}^-\cdot y_k$, one can get that $f=-2b-2$. Then we obtain
a representation $RA_{a,b}$ of ${\cal L}$ with basis
$\{x_i,y_i\,|\,i\in{\ma Z}\}$ and the actions: \vskip 0.2 cm
 $\!\!\!\!\!\!\! RA_{a,b}$\vspace*{-25pt}: \BE\label{a50}\hspace{1.2cm} \begin{array}{llllll} L_ix_j &\!\!\!=\!\!\!&
(a-j+ib)x_{i+j}, \!\!\!& L_iy_j &\!\!\!=\!\!\!& (a-j+i(b+\frac{1}{2}))y_{i+j}, \\
[0.2 cm] H_ix_j &\!\!\!=\!\!\!& -(2b+2)x_{i+j}, & H_iy_j
&\!\!\!=\!\!\!& -(2b+1)y_{i+j},
\\ [0.2 cm]
 G_i^-x_j &\!\!\!\!=\!\!\!& G_i^+y_j=0,  &  G_i^+x_j &\!\!\!=\!\!\!& \ y_{i+j},
 \ \  
 G_i^-y_j =
 2(a+i-j+2ib)x_{i+j},\!\!\!\!\!\!
\end{array}\EE (together with $cx_i=cy_i=0$ for all $i, j\in{\ma Z}$).

Obviously:
\\[0.2 cm] (i)\quad
 As $Vir$-modules, $V_0\cong A_{a, b}$ and $V_1\cong
A_{a, b'}$, where $b$ and $b'$ have some relations. \\[0.2 cm]
(ii)\quad For all $i\in{\ma Z}$, $H_i$ acts as constants on $V_0$
and $V_1$.
 \\ [0.3 cm]
{\bf Case 6.}\quad $a-b\in{\ma Z}$, then $a-b'\not\in{\ma Z}$.
Similar to the discussion in case 5, we have the following subcases:
\BC\begin{array}{lllll} \mbox{Subcase 6.1.} & f_{ij}=0, & f_{ij}'=1;
& b=0, & b'=\frac{1}{2}. \\ [0.3 cm] \mbox{Subcase 6.2.} & f_{ij}=0,
&
f_{ij}'=1; & b=-1, & b'=-\frac{1}{2}. \\
[0.3 cm] \mbox{Subcase 6.3.} & f_{ij}=f\not=0,  & f_{ij}'=f+1; & b=-\frac{1}{2}, & b'=0. \\
[0.3 cm] \mbox{Subcase 6.4.} & f_{ij}=f\not=0, &
f_{ij}'=\frac{a-i}{a-i-j}(f+1); & b=-\frac{1}{2}, & b'=0. \\[4pt] \mbox{Subcase 6.5.}
& f_{ij}=\frac{a-j}{a-i-j}f, & f_{ij}'=f+1; & b=0, & b'=\frac{1}{2}.
\\ [0.3 cm] \mbox{Subcase 6.6.} & f_{ij}=\frac{a-i-j}{a-j}f, &
f_{ij}'=f+1; &  b=-1, & b'=-\frac{1}{2}. \\
[0.3 cm] \mbox{Subcase 6.7.} & f_{ij}=f, &
f_{ij}'=\frac{a-i-j}{a-j}(f+1); & b=-\frac{3}{2}, &
b'=-1.\end{array}\EC
 Again by (\ref{a49}), we know that only Subcases 6.2 and 6.3 can
 occur. It is not difficult to see that they are contained in (\ref{a50}).
\\ [0.2 cm]
{\bf Case 7.}\quad $a-b'\in{\ma Z}$, then $a-b\not\in{\ma Z}$.
Similar to the discussion of case 6,  we obtain that if $a_{ij}^+=
1, \ b_{ij}^- =
 2(a-j+2i(b+\frac{1}{2})),\ a_{ij}^-=b_{ij}^+=0,\
 b'=b+\frac{1}{2}$, the module $V$ has the form of $RA_{a,b}$, for some $a, b\in{\ma C}  $.

 Similarly, we can write the other three cases of (\ref{a47}): \vskip 0.1 cm
 $\!\!\!\!\!\! RB_{a,b}$:\vspace*{-25pt}\  \BE\label{a71} \hspace{1.2cm} \begin{array}{llllll}L_ix_j &\!\!\!=\!\!\!&
(a-j+ib)x_{i+j}, \!\!\!\!& L_iy_j &\!\!\!=\!\!\!& (a-j+i(b+\frac{1}{2}))y_{i+j}, \\
[0.2 cm] H_ix_j &\!\!\!=\!\!\!& (2b+2)x_{i+j}, \!\!\!\!& H_iy_j
&\!\!\!=\!\!\!& (2b+1)y_{i+j},
\\ [0.2 cm]
 G_i^+ x_j \!\!&\!\!\!=\!\!\!& G_i^- y_j=0, \!\!\!\!&  G_i^-x_j &\!\!\!=\!\!\!& \ y_{i+j},
 \ \ 
 G_i^+y_j\!=
 2(a+i-j+2ib)x_{i+j},\!\!\!
 \end{array}\EE
when $a_{ij}^-= 1, \ b_{ij}^+ =
 2(a-j+2i(b+\frac{1}{2})),\ a_{ij}^+=b_{ij}^-=0,\
 b'=b+\frac{1}{2}$.
 \vskip 0.2 cm
 $\!\!\!\!\!\!\! RA'_{a,b}$:\vspace*{-27pt}\quad   \BE\label{a72} \hspace{1cm} \begin{array}{lllllll} L_ix_j &\!\!\!=\!\!\!&
(a-j+ib)x_{i+j},\!\!\!\!& L_iy_j &\!\!\!=\!\!\!& (a-j+i(b-\frac{1}{2}))y_{i+j}, \\
[0.2 cm] H_ix_j &\!\!\!=\!\!\!& -2b\ x_{i+j},  & H_iy_j
&\!\!\!=\!\!\!& -(2b+1)y_{i+j},
\\ [0.2 cm]
 G_i^+ x_j &\!\!\!=\!\!\!& G_i^- y_j=0, &   G_i^+y_j &\!\!\!=\!\!\!& \ x_{i+j},
 \ \ \ 
 G_i^-x_j=
 2(a-j+2ib)y_{i+j},\!\!\!\!\!\!
\end{array}\EE
when $a_{ij}^- =
 2(a-j+2ib), \ b_{ij}^+= 1, \  a_{ij}^+=b_{ij}^-=0,\
 b'=b-\frac{1}{2}$.
 \vskip 0.2 cm
 $\!\!\!\!\!\!\! RB'_{a,b}$:\vspace*{-27pt} \quad \BE\label{a73}\hspace{1cm}
  \begin{array}{llllll} L_ix_j &\!\!\!=\!\!\!&
(a-j+ib)x_{i+j}, \!\!\!& L_iy_j &\!\!\!=\!\!\!& (a-j+i(b-\frac{1}{2}))y_{i+j}, \\
[0.2 cm] H_ix_j &\!\!\!=\!\!\!& 2b\ x_{i+j}, & H_iy_j
&\!\!\!=\!\!\!& (2b+1)y_{i+j},
\\ [0.2 cm]
 G_i^- x_j &\!\!\!=\!\!\!& G_i^+ y_j=0, &  G_i^-y_j &\!\!\!=\!\!\!& \ x_{i+j},
 \ \ 
 G_i^+x_j =
 2(a-j+2ib)y_{i+j},
\end{array}\EE
when $a_{ij}^+ =
 2(a-j+2ib),\ b_{ij}^-= 1, \  a_{ij}^-=b_{ij}^+=0,\
 b'=b-\frac{1}{2}$.

 It is not
difficult to see that $RA_{a,b}\cong RA'_{a,b+\frac{1}{2}}$,
$RB_{a,b}\cong RB'_{a,b+\frac{1}{2}}$.
\\ [5pt]
{\bf \S3.2}\quad Now we consider all the possible deformations of
the representations which defined in section 3.1.
\\ [4pt]
{\bf Case 1.}\quad
 Suppose $V$ is an
indecomposable module which has the same composition factors as
those of $RA_{a,b}$ (in this case $RA_{a,b}$ is reducible). Let $V'$
be a non-zero submodule of $V$.
\\ [0.3 cm] {\bf Subcase 1.1.}\quad There exists $x_i\in V'$ for some
$i\in{\ma Z}$. By the sixth equation of (\ref{a50}), we obtain that
$y_j\in V'$ for all $j\in{\ma Z}$. By the last equation of
(\ref{a50}), \BE\label{a80}\begin{array}{lll} G_i^-y_{k-i}
&\!\!\!=\!\!\!&
2(a+i-k+i+2ib)x_{k}
=
2(a-k+2(b+1)i)x_k.\end{array} \EE  we see that $V'$
is a proper submodule of
 $V$ if and only if $a=k_0, \ b=-1$ for some $k_0\in{\ma Z}$.
In this case we can suppose $a=0, \ b=-1$. Then
$$V'=\mbox{span}_{\ma C}
 \{x_i, y_j\mid i, j\in{\ma Z},\ i\not=0 \}$$ is a nontrivial irreducible submodule of $V$,
 with the following relations:
\BE\label{a81}\begin{array} {llllll} L_ix_k &\!\!\!=\!\!\!&
-(i+k)x_{i+k}, & L_iy_j &\!\!\!=\!\!\!& -(\mbox{$\frac{i}{2}$}+j)y_{i+j}, \\
[0.2 cm] H_ix_k &\!\!\!=\!\!\!& 0, & H_iy_j &\!\!\!=\!\!\!& y_{i+j},
\\ [0.2 cm]  G_i^-x_k &\!\!\!=\!\!\!& G_i^+y_j=0, &  G_i^+x_k &\!\!\!=\!\!\!& y_{i+k},
\ \ \ \ 
G_i^-y_j =
-2(i+j)x_{i+j}, \end{array}\EE for all $j,
k\in{\ma Z}$ and $k\not=0$.

 In order to determine all possible
actions on $V$, we suppose that
$$  L_ix_0=l_ix_i, \ \ H_ix_0=h_ix_i, \ G_i^\pm x_0=g_i^{\pm}y_i \ \
\mbox{ for all}\ \ i\in{\ma Z}.$$ Act $[ L_{i-1}, H_1 ]=-H_i$ on
$x_0$, we can obtain that \BE\label{h} h_i=ih_1 ~~ \mbox{for all} ~~
i\in{\ma Z}. \EE Applying $[ L_i, G_j^- ]= (\frac{i}{2}-j)G^-_{i+j}$
to $x_0$, we obtain that \BE\label{g-}
-(\mbox{$\frac{i}{2}$}+j)g^-_j=(\mbox{$\frac{i}{2}$}-j)g^-_{i+j}.
\EE Set $j=0$, we get that $-\frac{i}{2}g^-_0=\frac{i}{2}g_i^-$,
i.e., $ g^-_i=-g^-_0 ~ \mbox{for} ~ i\in{\ma Z} \setminus \{0\}. $
Set $j=-i$ in (\ref{g-}), we get that
$-(\frac{i}{2}-i)g_{-i}^-=(\mbox{$\frac{i}{2}$}+i)g_0^-$, i.e., $
g_{-i}^-=\frac{3}{2}g_0^- ~ \mbox{for} ~ i\in{\ma Z} \setminus
\{0\}. $ Then we must have that \BE\label{g} g_i^-=0 ~~ \mbox{for
all} ~~ i\in{\ma Z}.\EE Applying $[ L_i, G_j^+ ]=
(\mbox{$\frac{i}{2}$}-j)G^+_{i+j}$ to $x_0$, we obtain that
\BE\label{l_i}
l_i=-\mbox{$(\frac{i}{2}+j)g_{j}^+-(\frac{i}{2}-j)g^+_{i+j}.$} \EE
Following $[ H_i, G^+_j ]\cdot x_0= G^+_{i+j}\cdot x_0$, we get that
\BE\label{a84} h_i=g_{j}^+-g_{i+j}^+ ~~ \mbox{for all} ~~ i,
j\in{\ma Z}. \EE Set $j=0$, and by (\ref{h}), one can get that
\BE\label{b1}g^+_i=g^+_0-ih_1.\EE It follows (\ref{l_i}) that
\BE\label{l}\mbox{$ l_i=\frac{i^2}{2}h_1-ig^+_0.$}\EE If $h_1=0$,
then
$$l_i=-ig^+_0, ~~ h_i=0,~~g^-_i=0, ~~ g^+_i=g^+_0,$$ and it satisfies
(\ref{a81}) (rescaling $x_0$ by $(g^+_0)^{-1}x_0$). Then it is not a
deformation of $V$. Hence we suppose that $h_1\not=0$. Rescaling
$x_0$ by $h_1^{-1}x_0$, and follows (\ref{l_i}), (\ref{b1}), we can
obtain that \BE\label{b2}\mbox{$ l_i =
-i\frac{g^+_0}{h_1}+\frac{i^2}{2}, ~~~ h_i = i, ~~~ g^+_i
=\frac{g^+_0}{h_1}-i, ~~~ g^-_i = 0.$} \EE Then we get a deformation
of $RA_{a, b}$, denoted by $RA_\alpha$, which is an indecomposable
module with the following relations: \vskip 0.2 cm $
RA_\alpha$\vspace* {-25pt}: \BE\label{b3}\ \ \ \ \ \ \hspace*{20pt}
\begin{array}{llllllll}
L_ix_k=-(i+k)x_{i+k}, & L_ix_0= (i\alpha+\frac{i^2}{2})x_i,
 & L_iy_j= -(\frac{i}{2}+j)y_{i+j},\\[4pt]
H_iy_j=y_{i+j},&H_ix_k= 0,  &  H_ix_0=i x_i,  \\[4pt]
 G^-_ix_j=G^+_iy_j=0, &G^-_iy_j= -2(i+j)x_{i+j},
 \\[4pt] G^+_ix_k=y_{i+k}, &  G^+_ix_0=-(\alpha+i)y_i,
\end{array}\EE
where $\alpha=-\frac{g_0^+}{h_1}\in{\ma C}$, $j, k\in{\ma Z}$ and
$k\not=0$.
\\ [0.2 cm] {\bf Subcase 1.2.}\quad There exists $y_{j_0}\in V'$ for
some $j_0\in{\ma Z}$. By (\ref{a50}), in order to make $V'$ is a
proper submodule of $V$, we must have that for all $i\in{\ma Z}$,
\BC G^-_iy_{j_0}=2(a+i-j_0+2ib)x_{i+j_0}=0, \EC it follows that
$a=j_0, b=-\frac{1}{2}$. Without loss of generality, we can suppose
that $j_0=0$, then we have $a=0, b=-\frac{1}{2}$. Therefore,
$V'={\ma C}y_0$ is a trivial proper submodule of $V$, and set
$V''=$span$_{\ma C}\{ x_i, y_j\mid i, j\in{\ma Z}, j\not=0 \}$. We
have the following relations: \BE\label{b5}& &
\begin{array}{llllll} L_ix_j &\!\!\!=\!\!\!&
-(\frac{i}{2}+j)x_{i+j}, & L_iy_k &\!\!\!=\!\!\!& -ky_{i+k}, \\
[0.2 cm] H_ix_j &\!\!\!=\!\!\!& -x_{i+j}, & H_iy_k &\!\!\!=\!\!\!&
0,\\ [0.2 cm]
G_i^-x_k &\!\!\!=\!\!\!& G_i^+y_j=0, & G_i^+x_k &\!\!\!=\!\!\!& \ y_{i+k},
\ \ \ 
G_i^-y_j=
-2jx_{i+j},
\end{array}\EE
where $j, k\in{\ma Z}$ and $k\not=-i$. Suppose that \BC
L_iy_{-i}=l_iy_0, ~~~ H_iy_{-i}=h_iy_0, ~~~
G^\pm_ix_{-i}=g^\pm_iy_0.\EC  Similar to the arguments about Case 1,
one can get that $$ l_i=ig^+_0+\mbox{$\frac{i^2}{2}$}h_1, ~~
g^+_i=g^+_0+ih_1, ~~ g^-_i=0, ~~ h_i=ih_1.$$ If $h_1=0$, then
$$l_i=ig^+_0, ~~ g^+_i=g^+_0, ~~ g^-_i=0, ~~
h_i=0.$$ It is not a deformation too, so we suppose that
$h_1\not=0$.
 Rescaling $y_0$ by
$h_1y_0$, we can get a new representation of ${\cal L}$, we denote
it by $RA^\beta$.
 \vskip 0.2 cm
 $RA^\beta$:\quad  satisfies
(\ref{b5}) and the following relations (set
$\beta=\frac{g_0^+}{h_1}$): \BE\label{b7}
L_iy_{-i}=(i\beta+\mbox{$\frac{i^2}{2}$})y_0, ~~~ H_iy_{-i}=iy_0,
~~~ G^-_ix_{-i}=0, ~~~ G^+_ix_{-i}=(\beta+i)y_0. \EE \vskip 0.1 cm
Obviously, $RA_\alpha$ has a nontrivial submodule with codimention
one, and $RA^\beta$ has a trivial submodule with dimention one.
\\ [0.2 cm]   {\bf Case 2.}\quad Now we discuss the
deformations of $RB_{a,b}$. Since the discussion is similar to Case
1, we will not give the detail, and only enumerate the results.
$RB_{a,b}$ also has two deformations: \vskip 0.2 cm $
\!\!\!RB_\alpha$\vspace* {-25pt}: \BE\label{c1}\ \ \ \ \ \ \ \ \ \ \
\
\begin{array}{llllllll}
 L_ix_k=-(i+k)x_{i+k}, &L_ix_0=
-(i\alpha+\frac{i^2}{2})x_i,
 & L_iy_j= -(\frac{i}{2}+j)y_{i+j},\\[4pt]
 H_iy_j=-y_{i+j},
 & H_ix_k= 0,  &  H_ix_0= i x_i,  \\ [4pt]
 G^+_ix_j=G^-_iy_j=0, & G^+_iy_j=-2(i+j)x_{i+j}, \\[0.2 cm]
  G^-_ix_k=y_{i+k}, &  G^-_ix_0=(\alpha+i)y_i,
\end{array}\EE
where $\alpha=\frac{g_0^-}{h_1}\in{\ma C}$, $h_1\not=0$, $j,
k\in{\ma Z}$, and $k\not=0$.  \vskip 0.2 cm
 $RB^\beta$\vspace*{-24pt}: \BE\label{c3}\ \ \ \ \ \ \ \ \ \ \  \
\begin{array}{llllllll} L_ix_j =-(\frac{i}{2}+j)x_{i+j},  & L_iy_k= -ky_{i+k}, &
L_iy_{-i}= -(i\beta+\frac{i^2}{2})y_0, \\[4pt]
H_ix_j= x_{i+j},&  H_iy_k= 0,  & H_iy_{-i}=iy_0, \\
[0.2 cm]
G_i^+x_j=G_i^-y_j=0, &  G_i^-x_k= y_{i+k},\\
[0.2 cm]  G_i^-x_{-i}= -(\beta+i)y_0, &  G_i^+y_j = -2jx_{i+j},
\end{array}\EE
where $\beta=-\frac{g_0^-}{h_1}\in{\ma C}$, $h_1\not=0$, $j,
k\in{\ma Z}$ and $k\not=-i$. \vskip0.3cm

This completes the proof of Theorem 1.1.


\begin{thebibliography}{9999}\parskip0pt\lineskip0pt\small
  \bibitem {A1} M. Ademollo, L. Brink, A. d'Adda, R. Auria,
 E. Napolitano, S. Sciuto, E. del Giudice, P. di Vecchia,
 S. Ferrara, F. Gliozzi, R. Musto,  R. Pettorino,
  Supersymmetric strings and colour confinement. Phys. Lett. B 62,
  105(1976).
  \bibitem {A} E. Arbarello, C. De Concini, V.G. Kac and C. Procesi,
  Moduli spaces of curves and representation theory, Comm. Math.
  Phys., 117(1), 1-36(1988).
     \bibitem {CK} S.L. Cheng,  V.G. Kac, A new $N=6$ superconformal
 algebra. Comm. Math. Phys. 186, 219-231(1997).
    \bibitem {K1} V.G. Kac, Lie superalgebras. Adv. Math. 26, 8-97(1977).
  \bibitem {KL} V.G. Kac,  J.W. van de Leuer, {\it On classification of superconformal
 algebras.} Strings 88, Sinapore: World Scientific, (1988).
  \bibitem {IL} I. Kaplansky, L.J. Santharoubane, Harish-Chandra modules over the Virasoro algebras. in MSRI
 Publ.4, 217-231(1987).
  \bibitem {KE} E. Kiritsis, Character formula and the structure of
 the represetations of the $N=1$, $N=2$ superconformal algebrass.
 Int. J. Mod. Phys. A 3, 1871-1906(1988).
\bibitem {LZ} R. Lu, K. Zhao, Classification of Irreducible Weight Modules over the Twisted Heisenberg-Virasoro
 Algebras. Arxiv: math.ST/0510194 v1, 10 Oct 2005.
 \bibitem {CA} C. Martin, A. Piard, Indecomposable modules for the Virasoro Lie algebra and a conjecture of Kac. Comm. Math. Phys. 137,
 109-132(1991).
  \bibitem {S1} Y. Su, Classification of Harish-Chandra Modules over the Super-Virasoro Algebras. Comm. Alg. 23(10),
 3653-3675(1995).
 \bibitem {S2} Y. Su, A classification of indecompable $sl_2({\ma C})$-modules and a conjecture of Kac on irreducible modules over
 the Virasoro algebra. J. Algebra, 161(1), 33-46(1993).
 \bibitem {SZ} Y. Su,  K. Zhao, Generalized Virasoro and super-Virasoro algebras and modules of the intermediate series.
  J. Algebra,
 252(1), 1-19(2002).
\end{thebibliography}
    \end{document}